\documentclass[12pt]{article}
\usepackage[all]{xy}
\usepackage{amsmath}
\usepackage{amsfonts}
\usepackage{amssymb}
\usepackage{amscd}
\usepackage{amsthm}
\usepackage{latexsym}
\usepackage{amsbsy}
\newtheorem{lem}{Lemma}[section]
\newtheorem{prop}{Proposition}[section]

\newtheorem{definition}{Definition}[section]
\newtheorem{example}{Example}[section]

\newcommand{\ot}{\otimes}

\begin{document}
\title{Introduction to  Hopf-Cyclic Cohomology}
\author {Masoud Khalkhali and Bahram Rangipour}
\date{}
\maketitle
\begin{abstract} We review the recent progress in the study of cyclic
cohomology in the presence of  Hopf symmetry.

 \end{abstract}

\section{Introduction}
In their study of the index theory of transversally elliptic operators \cite{cm2},  Connes and Moscovici
 develope a cyclic cohomology theory for Hopf algebras which
 can be regarded, {\it post factum},  as the right  noncommutative  analogue    of group homology
  and Lie algebra homology.   One of the main reasons
  for  introducing this theory was to obtain a
  {\it noncommutative characteristic map}
  $$\chi_{\tau}: {HC}_{(\delta,\sigma)}^\ast (H) \longrightarrow  HC^\ast (A),$$
  for an action of a  Hopf
  algebra $H$ on an  algebra $A$ endowed with an ``invariant trace'' $\tau:
  A \to \mathbb{C}$. Here, the pair $(\delta, \sigma)$ consists of a grouplike element $\sigma \in H$ and
   a characater $\delta: H\to \mathbb{C}$ satisfying certain compatibility conditions explained in
   Section 2.3 below.

In \cite{kr2} we found a new approach to this subject and extended it, among other things, to
a cyclic cohomology theory  for triples $(C, H, M)$, where $C$ is  a coalgebra endowed with an action of
a Hopf algebra $H$ and  $M$ is an $H$-module and an $H$-comodule  satisfying
some extra compatibility conditions.
It was observed that
the theory  of Connes and Moscovici corresponds to $C=H$ equipped with
the regular action of $H$ and $M$ a one dimensional $H$-module with an extra
structure.

 One of the main ideas of $\cite{kr2}$ was to view the
Hopf-cyclic cohomology as the  cohomology of the {\it invariant} part of
certain natural complexes attached
to $(C, H, M)$. This is  remarkably similar to interpreting the cohomology
 of the Lie algebra of a
 Lie group as the invariant part of the de Rham cohomology of the Lie group. The second main
idea was to introduce {\it coefficients} into the theory.   This also explained the
important role played by the
so called modular pair $(\delta, \sigma)$ in \cite{cm2}.

Since the module $M$ is a noncommutative analogue of coefficients
 for  Lie algebra and group homology theories, it is of utmost importance
 to understand  the most general type of coefficients allowable. In fact the periodicity
condition $\tau_n^{n+1}=id$ for the cyclic operator and the fact that all simplicial and
cyclic operators have to descend to the invrainat complexes, puts very
stringent conditions on the type of the $H$-module $M$. This problem
that remained
 unsettled in our paper \cite{kr2}  is completely solved
 in Hajac-Khalkhali-Rangipour-Sommerh\"auser papers \cite{hkrs1, hkrs2} by
 introducing the class of {\it stable anti-Yetter-Drinfeld modules} over a
 Hopf algebra. The category of anti-Yetter-Drinfeld modules over a Hopf algebra $H$ is a
 twisting, or `mirror image' of the category of Yetter-Drinfeld $H$-modules.
 Technically it is obtained from the latter  by replacing the antipode $S$ by $S^{-1}$ although this
 connection is hardly illuminating.

 In an effort to make this paper more accessible,  we cover basic background material,
 with many examples,  on
 Hopf algebras and the emerging role of {\it Hopf symmetry} in noncommutative geometry and its applications
 \cite{cm2, cm3, cm4, cm5}. This is justified since certain doses of the
 {\it ``yoga of Hopf algebras''},
 in the noncommutative and non-cocommutative case, is necessary to understand these works. Following these works, we emphasize the universal role played by the
 Connes-Moscovici Hopf algebra $\mathcal{H}_1$ and its cyclic cohomology
 in applications of noncommutative geometry to transverse geometry and
 number theory. See also Marcolli's article \cite{mar} as well as Connes-Marcolli
 articles \cite{cma1, cma2}
  and references therein  for  recent 
  interactions between number theory
 and noncommutative geometry.

It is a great pleasure to thank Alain Connes  and Henri Moscovici
for their  support and many illuminating discussions, and Matilde
Marcolli for organizing the MPI, Bonn, conferences on noncommutative
geometry and number theory and for her support.
  We are also much
obliged to our collaborators  Piotr M.  Hajac and Yorck Sommerh\"auser.

 \section{Preliminaries}

In this section we recall some  aspects of Hopf algebra theory  that are most
relevant to the current status of Hopf-cyclic cohomology theory.

\subsection{Coalgebras, bialgebras,  and Hopf algebras}

We assume our Hopf algebras, coalgebras, and algebras are over a fixed  field $k$ of
characteristic zero. Most of our definitions and  constructions however continue to work
over an arbitrary commutative
ground ring $k$. Unadorned $\otimes$ and $Hom$  will always be over k and $I$ will always denote
an
identity map whose domain will be clear from the context.
By a {\it coalgebra} over $k$ we mean a $k$-linear space $C$ endowed
with $k$-linear maps
$$\Delta : C\longrightarrow C\otimes C, \quad \varepsilon :
C\longrightarrow k,$$
called {\it comultiplication} and {\it counit} respectively, such that
$\Delta $ is {\it coassociative} and $\varepsilon $ is the counit of
$\Delta$. That is,
$$(\Delta \otimes I)\circ \Delta =(I\otimes \Delta)\circ \Delta :
C\longrightarrow C\otimes C \otimes C,$$
$$ (\varepsilon \otimes I)\circ \Delta =(I\otimes \varepsilon)\circ \Delta
=I.$$
$C$ is called {\it cocommutative} if $\tau \Delta =\Delta$, where
$\tau: C\otimes C \rightarrow C\otimes C$ is the {\it fillip} $x\otimes
y \mapsto y\otimes x$.

We use  Sweedler-Heynemann's notation for
comultiplication, with summation suppressed, and write
$$ \Delta (c)= c^{(1)}\otimes c^{(2)}.$$
With this notation the axioms
for a coalgebra   read as
$$c^{(1)}\otimes c^{(2)(1)} \otimes c^{(2)(2)}=
c^{(1)(1)}\otimes c^{(1)(2)} \otimes c^{(2)},$$
$$\varepsilon
(c^{(1)})(c^{(2)})=c=(c^{(1)})\varepsilon(c^{(2)}),$$
for all $c\in C$.
We put
$$c^{(1)}\otimes c^{(2)}\otimes c^{(3)}:=(\Delta \otimes I)\Delta (c).$$
Similarly, for  {\it iterated comultiplication} maps
$$\Delta ^n: =(\Delta \otimes I)\circ \Delta^{n-1}: C\longrightarrow C^{\otimes (n+1)},\quad \Delta^1=\Delta,$$
we write
$$\Delta^n (c)= c^{(1)}\otimes \cdots \otimes c^{(n+1)},$$
where summation is understood. Many notions for  algebras have
their dual analogues for coalgebras. Thus, one can easily define
such notions  like, {\it subcoalgebra}, (left, right, two sided)
{\it coideal}, {\it quotient coalgebra}, and {\it morphism of
coalgebras} \cite{maj, sw}.

A left $C$-{\it comodule} is a linear space $M$ together with a linear
map $\rho: M\to C\otimes M$ such that $(\rho \otimes 1)\rho =\Delta
\rho$ and $(\varepsilon \otimes 1)\rho =\rho$. We write
$$ \rho(m)=m^{(-1)} \otimes
m^{(0)},$$
where summation is understood,
to denote the coaction $\rho$. Similarly if $M$ is a right $C$-comodule, we
write
$$\rho (m)= m^{(0)} \otimes m^{(1)}$$
to denote its coaction
$\rho: M \rightarrow M\otimes C$. Note that module elements are always assigned zero index. Let
$M$ and $N$ be left $C$-comodules. A $C$-{\it colinear map} is
a linear map $f: M \to N$ such that $\rho_N f= (1\otimes f)\ \rho_M$. The category of left $C$-comodules
is an abelian category; note that, unlike the situation for algebras,
for this to be true, it is important that $C$ be a flat $k$-module
which will be the case if $k$ is a field.

Let $C$ be a coalgebra, $A$ be a unital  algebra, and $ f , g :C
\rightarrow A$ be $k$-linear maps. The {\it convolution product}
of $f$ and $g$, denoted by  $f*g$, is defined as the composition
$$ C\overset{\Delta}{\longrightarrow} C\otimes C \overset{f\otimes g}{\longrightarrow} A\otimes A,$$
or equivalently  by
$$(f*g)(c)=f(c^{(1)})g(c^{(2)}).$$
It is easily checked that under the convolution product
$Hom(C, A)$ is an associative unital algebra.
 Its unit  is the map $
e :C \rightarrow A$, $e(c)= \varepsilon (c)1_A$.
In particular the linear dual of a coalgebra $C$, $C^*= Hom (C, k),$ is an
algebra.

A {\it bialgebra} is a unital  algebra endowed with a compatible coalgebra
structure. This means that the coalgebra structure maps $\Delta : B\longrightarrow B\otimes B
, \quad \varepsilon :
B\longrightarrow k,$ are morphisms of unital algebras. This is equivalent to multiplication
and unit maps of $A$  being morphisms of coalgebras.

A {\it Hopf algebra} is a bialgebra endowed with an {\it antipode}. By definition, an antipode
for  a
 bialgebra $H$
is a linear map $S :H \rightarrow H$ such that
$$S*I=I*S=\eta \varepsilon,$$
where $\eta: k \to H$ is the unit map of $H$. Equivalently,
$$  S(h^{(1)}) h^{(2)}= h^{(1)}S(h^{(2)})= \varepsilon (h) 1,$$
 for all $h \in H$.
Thus $S$ is the inverse of the identity map $I: H \to H$ in the convolution
algebra $Hom (H, H)$. This shows that the antipode is unique, if it
exists at all.  The following properties of the antipode are well
known:

1. If $H$ is commutative or cocommutative then $S^2=I$. The converse need
not be true.\\
2. S is an anti-algebra map and  an anti-coalgebra map. The latter
means
$$ S(h^{(2)}) \otimes S(h^{(1)})=S(h)^{(1)} \otimes S(h)^{(2)},$$
for all $h\in H$.

We give a few examples of Hopf algebras:\\

1. Commutative or  cocommutative Hopf algebras are closely related to
groups and Lie algebras. We give a few examples to indicate this
connection\\

1.a. Let $G$ be a discrete group (need not be finite) and $H=kG$
the group algebra of $G$ over $k$. Let
$$ \Delta (g)=g\otimes g, \quad S(g)=g^{-1}, \quad \text{and}\;
\varepsilon (g)=1,$$ for all $g\in G$ and linearly extend them to
$H$. Then it is easy to check that $(H, \Delta, \varepsilon, S)$
is a cocommutative Hopf
algebra. It is of course commutative if and only if $G$ is commutative.\\

1.b. Let $\mathfrak{g}$  be a $k$-Lie algebra and
$H=U(\mathfrak{g})$ be  the universal enveloping algebra of
$\mathfrak{g}$. Using the universal property of $U(\mathfrak{g})$
one checks that there are uniquely defined  algebra homomorphisms
$\Delta: U(\mathfrak{g})\to U(\mathfrak{g})\otimes
U(\mathfrak{g})$, $\varepsilon: U(\mathfrak{g})\to k$ and an
anti-algebra map $S:U(\mathfrak{g})\to U(\mathfrak{g})$,
determined by
$$ \Delta (X)=X\otimes 1 +1\otimes X, \quad \varepsilon (X)=0, \quad \text{and} \; \; S(X)=-X,$$
for all $X \in \mathfrak{g}$.  One then checks easily that
$(U(\mathfrak{g}), \Delta, \varepsilon, S)$ is a cocommutative Hopf
algebra. It is commutative if and only if $\mathfrak{g}$ is an abelian Lie algebra,
in which case $U(\mathfrak{g})=S(\mathfrak{g})$ is the symmetric algebra of $\mathfrak{g}$. \\

1.c. Let $ G$ be a compact topological group. A continuous
function $f:G \to \mathbb{C}$ is called {\it representable} if the
set of left translates of $f$ by all elements of $G$ forms a finite
dimensional subspace of the space $C(G)$ of all continuous complex
valued functions on $G$. It is then easily checked that the set of
all representable functions,
 $H =Rep(G)$, is a subalgebra
 of the algebra of continuous functions on $G$. Let $m: G\times G \to G$ denote the multiplication of $G$
 and $m^*: C(G\times G)\to C(G), \; m^*f (x, y)=f(xy), $ denote its  dual map. It is easy to see that
  if $f$
 is representable, then
 $m^*f \in Rep(G)\otimes Rep (G) \subset C(G\times G)$.  Let $e$ denote  the identity of $G$.
 One can easily check
 that  the relations
 $$ \Delta f =m^*f, \quad \varepsilon f=f(e), \quad \text{and} \; \;
 (Sf)(g)=f(g^{-1}),$$
 define a Hopf algebra structure on $Rep (G)$. Alternatively, one can describe $Rep (G)$ as the linear
 span of matrix coefficients of all finite dimensional complex representations of $G$.
 By {\it Peter-Weyl's Theorem,} $Rep (G)$ is a dense
 subalgebra of $C(G)$. This algebra is finitely
 generated (as an algebra) if and only if $G$ is a Lie group.\\

1.d. The coordinate ring of an affine algebraic group $H=k[G]$ is
a commutative Hopf algebra. The maps $\Delta$, $\varepsilon$, and
$S$ are the duals of the multiplication, unit, and inversion  maps
of $G$, respectively. More generally, an {\it affine group
scheme}, over a commutative ring $k$, is a commutative Hopf
algebra over $k$. Given such a Hopf algebra $H$, it is easy to see
that for any commutative $k$-algebra $A$, the set $Hom_{Alg}(H,
A)$ is a group under convolution product and $A\mapsto
Hom_{Alg}(H, A)$ is a functor from the category $ComAlg_k$ of
commutative $k$-algebras to the category of groups. Conversely,
any representable functor $ComAlg_k \to Groups$ is represented by
a, unique up to isomorphism,  commutative $k$-Hopf algebra. Thus
the category  of affine group schemes is equivalent to the
category of representable
functors $ComAlg_k \to Groups.$  \\

2. Compact quantum groups and quantized enveloping algebras are
examples of noncommutative and noncommutative Hopf algebras
\cite{maj}. We won't recall these examples here. A very important
example for noncommutative geometry and its applications to
transverse geometry and number theory is the {\it Connes-Moscovici
Hopf algebra} $\mathcal{H}_1$ \cite{cm2, cm4, cm5} which  we
recall now.  Let $\mathfrak{g}_{aff}$ be the Lie algebra of the
group of affine transformations of the line with linear basis $X$
and $Y$ and the relation $ [Y,X] = X$. Let $\mathfrak{g}$ be an
abelian Lie algebra with basis $\{\delta_n; \; n=1, 2, \cdots\}$.
It is easily seen that $\mathfrak{g}_{aff}$ acts on $
\mathfrak{g}$ via
$$[Y, \delta_n] = n \delta_n,  \quad  [X, \delta_n] = \delta_{n+1},$$
for all $n$. Let $\mathfrak{g}_{CM}:=\mathfrak{g}_{aff}\rtimes \mathfrak{g}$ be the
corresponding semidirect product Lie algebra.
As an algebra, $\mathcal{H}_1$ coincides with
the universal enveloping algebra of the Lie algebra $\mathfrak{g}_{CM}$. Thus $\mathcal{H}_1$
is the universal algebra generated by  $\{X, Y, \delta_n ; n =
1, 2, \cdots \}$ subject to relations
$$[Y,X] = X, \; \;[Y, \delta_n] = n \delta_n, \; \;[X, \delta_n] =
\delta_{n+1},\; \;
[\delta_k, \delta_l] = 0,$$
 for $ n, k, l = 1 ,
2, \cdots.$
We let the counit of $\mathcal{H}_1$ coincide with the counit of $U(\mathfrak{g}_{CM})$. Its coproduct
and antipode, however, will be certain deformations of the coproduct and antipode of
$U(\mathfrak{g}_{CM})$ as follows.  Using the universal property of  $U(\mathfrak{g}_{CM})$,
one checks that
the relations
  $$\Delta Y = Y \otimes 1 + 1 \otimes  Y, \quad  \Delta \delta_1 =
  \delta_1 \otimes  1 + 1 \otimes  \delta_1, $$
$$\Delta X = X \otimes  1 + 1 \otimes  X + \delta_1 \otimes  Y,$$
determine a unique  algebra map $ \Delta : \mathcal{H}_1 \to
\mathcal{H}_1 \otimes  \mathcal{H}_1$. Note that $\Delta$ is not
cocommutative and it differs from the corrodent of the enveloping
algebra $U(\mathfrak{g}_{CM})$. Similarly, one checks that there
is a unique antialgebra map $S:  \mathcal{H}_1 \to \mathcal{H}_1$
determined by the relations
$$ S(Y ) = -Y, \; \;  S(X) = -X + \delta_1Y, \; \; S(\delta_1) = -\delta_1.$$
Again we note that this antipode also differs from the antipode of
$U(\mathfrak{g}_{CM})$, and in particular $S^2\neq I$. In fact
$S^n\neq I$ for all $n$. In the next section we will show,
following Connes-Moscovici \cite{cm2}, that $\mathcal{H}_1$ is a
bicrossed product of  Hopf algebras $U(\mathfrak{g}_{aff})$ and
$U(\mathfrak{g})^*$, where $\mathfrak{g}$ is a pro-nilpotent Lie
algebra to be described precisely in the next section.

 Let $H$ be a Hopf algebra. A {\it grouplike} element of $H$  is a non-zero element $g\in H$
 such that $\Delta g = g \otimes g$. We have, from the axioms for the antipode, $g S(g)=S(g)g=1_H$ which
 shows that $g$ is invertible. It is easily seen that
  grouplike elements of $H$ form a subgroup of the multiplicative group
  of $H$. For example, for $H=kG$ the set of grouplike elements coincide with the group $G$ itself.
  A {\it primitive element} is an element $x \in H$ such that $\Delta x
  =1\otimes x +x\otimes 1$. It is easily seen that the bracket $[x, y]:=xy-yx$
  of two primitive elements is again a primitive element. It follows that
  primitive elements form a Lie algebra. Using the {\it Poincare-Birkhoff-Witt } theorem, one shows that
  the set of primitive elements of  $H= U(\mathfrak{g})$  coincide with the Lie algebra
  $\mathfrak{g}$.

  A {\it character} of a Hopf algebra is a unital algebra map $\delta : H
\rightarrow k$. For example the counit $\varepsilon : H\to k$ is a character. For a less trivial
example,
 let $G$ be a {\it non-unimodular} real Lie group and
$H=U(\mathfrak{g})$ the universal enveloping algebra of the Lie algebra of $G$. The
{\it modular function} of $G$, measuring the difference between the right and left Haar measure on $G$,
 is
 a smooth group
homomorphism $\Delta: G \to \mathbb{R}^+$. Its derivative at identity defines
a Lie algebra map $\delta: \mathfrak{g} \to \mathbb{R}$. We denote its natural extension by
$\delta : U(\mathfrak{g})\to \mathbb{R}$. It is obviously a character of $U(\mathfrak{g}).$
  For a concrete example, let $G=Aff (\mathbb{R})$ be the group of affine
transformations of the real line. It is a non-unimodular  group with
modular homomorphism given by
$$\Delta \left( \begin{matrix} a &  b\\ 0&1   \end{matrix}\right)= \vert a\vert.$$
The corresponding infinitesimal character on
$\mathfrak{g}_{aff}=Lie(G)$ is given by
$$\delta (Y)=1, \quad \delta (X)=0,$$ where $Y=\left( \begin{matrix} 1 & 0\\ 0&0  \end{matrix}\right)$
and
$X=\left( \begin{matrix} 0 & 1\\ 0&0  \end{matrix}\right)$ are a basis for
$\mathfrak{g}_{aff}.$ We will see that this character plays an important role
in constructing a {\it twisted antipode} for the Connes-Moscovici Hopf
algebra $\mathcal{H}_1$.

If $H$ is a Hopf algebra, by a left $H$-module (resp. left $H$-comodule), we mean a
left module (resp. left comodule) over the underlying algebra (resp. the underlying coalgebra)
of $H$.

Recall that a {\it monoidal}, or {\it tensor}, category $(\mathcal{C},
\otimes, U, a, l, r)$ consists of a category $\mathcal{C}$, a
functor $\otimes: \mathcal{C}\times \mathcal{C}\to \mathcal{C}$, an
object $U\in \mathcal{C}$ (called the {\it unit object}), and natural
isomorphisms
$$ a=a_{A, B, C}: A\otimes (B\otimes C) \to (A\otimes
B)\otimes C,$$
$$ l=l_A: U\otimes A \to A, \quad \quad r=r_A: A\otimes U \to A,$$
(called the {\it associativity} and {\it unit constraints}, respectively) such that the so
called {\it pentagon}
and {\it triangle} diagrams commute:\\

\begin{center}
$\begin{xy}
\xymatrix{ &((A\otimes B)\otimes C)\otimes D\ar[dl] \ar[dr]&\\
(A\otimes (B \otimes C))\otimes D\ar[d]& & (A\otimes B) \otimes (C \otimes D)\ar[d]\\
A\otimes ((B \otimes C)\otimes D )\ar[rr]& & A\otimes (B \otimes
(C\otimes D))}
\end{xy}$
\end{center}

\begin{center}
$\begin{xy}
\xymatrix{ (A\otimes U)\otimes B  \ar[rr] \ar[dr]& & A\otimes (U \otimes B)\ar[dl]\\
& A\otimes B& }
\end{xy}$
\end{center}

  The coherence theorem of MacLane \cite{ml}  guarantees
that all diagrams composed from $a, l, r$ by tensoring, substituting and composing, commute.\\

A {\it braided monoidal} category, is a monoidal category $\mathcal{C}$
endowed with a natural family of isomorphisms
$$c_{A, B}: A\otimes B \to
B\otimes A,$$ called  {\it braiding}  such that the following
diagram and the one obtained from it by
replacing $c$ by $c^{-1}$ commute ({\it Hexagon axioms}):\\

\begin{center}
$\begin{xy}
\xymatrix{ & A\otimes (B\otimes C)\ar[r] & (B\otimes C)\otimes A\ar[dr] &\\
(A\otimes B)\otimes C\ar[ur] \ar[dr]& & &(C\otimes B)\otimes A \\\
&(B\otimes A)\otimes C\ar[r]& C\otimes (B\otimes A) \ar[ru] &}
\end{xy}$
\end{center}

A braiding is called a {\it symmetry} if we have
$$c_{A, B}\circ c_{B, A}=I$$
for all $A, B$. A {\it symmetric monoidal category}, is a monoidal
category endowed with a symmetry.

Let $H$ be a bialgebra. Then thanks to the existence of a comultiplication on $H$, the category
$H-Mod$ of
left $H$-modules is a    monoidal category:
if $M$ and $N$ are left $H$-modules, their
 tensor product over $k$, $M\otimes N$, is again an $H$-module via
 $$h(m\otimes n)= h^{(1)}m\otimes h^{(2)}n  \label{tp}.$$
One can easily check that with associativity constraints defined as the
natural isomorphism of $H$-modules $M\otimes (N \otimes P) \to (M\otimes N)\otimes P$,
$m\otimes (n\otimes p)\mapsto (m\otimes n)\otimes p$,
and with unit object $U=k$ with trivial $H$-action via the counit $\varepsilon$,
$$h(1)=\varepsilon (h)1,$$
$H-Mod$ is  a monoidal category. If $H$ is cocommutative, then one checks that the map
$$c_{M, N}: M\otimes N \to N\otimes M, \quad  c_{M, N}(m\otimes
n)=n\otimes m,$$
is a morphism of $H$-modules and is a symmetry operator on $H-Mod$,
turning it into a symmetric monoidal
category.

The category $H-Mod$ is not braided in general. For that to happen,  one must either restrict the
class of modules to what is
called Yetter-Drinfeld modules, or restrict the class of Hopf algebras to
 quasisimilar Hopf algebras to obtain a braiding on $H-Mod$
 \cite{maj}. We will discuss the first scenario in the next section and
 will see that,
 quite unexpectedly, this question is closely related to Hopf-cyclic cohomology.

Similarly, the category $H-Comod$  of left $H$-comodules is a monoidal category:
if $M$ and $N$ are left $H$-comodules, their tensor product $M\otimes N$
is again an $H$-comodule
via $$\rho (m\otimes n)=m^{(-1)} n^{(-1)}\otimes m^{(0)}\otimes n
^{(0)}.$$ Its unit object is $U=k$ endowed with the $H$-coaction $r\in k
\mapsto 1_H\otimes r.$ If $H$ is commutative
then $H-Comod$ is a
symmetric monoidal category. More generally, when $H$ is co-quasitriangular,
$H-Comod$ can be endowed with a braiding \cite{maj}.

\subsection{Symmetry in noncommutative geometry}

The idea of {\it symmetry} in noncommutative geometry is encoded
by the {\it action} or {\it coaction}  of a Hopf algebra on an
algebra or on a coalgebra. Thus there are four possibilities in
general that will be referred to as (Hopf-) {\it module algebra},
{\it module coalgebra}, {\it comodule algebra}, and {\it comodule
coalgebra} (for each type, of course, we still have a choice of
left or right action or coaction). We call them {\it symmetries of
type A, B, C} and  {\it D}, respectively. We will see in the next
sections that associated to each type of symmetry there is a
corresponding cyclic cohomology theory with coefficients. These
theories in a certain sense are generalizations of equivariant de
Rham cohomology with coefficients in an equivariant local system.

 Let $H$ be a Hopf algebra. An algebra $A$ is called a left $H$-\textit{module algebra},
 if it  is a left $H$-module and the multiplication  map   $A\otimes A\to A$ and the unit map
 $k \to A$ are   morphisms of  $H$-modules. This means
$$ h(ab)=h^{(1)}(a)h^{(2)}(b), \quad \text{and}\quad h(1)=\varepsilon (h)1,$$
for all $h \in H$ and $a, b \in A$ (summation is understood).  Using the relations
$\Delta g=g\otimes g$ and $\Delta x=1\otimes x +x \otimes 1$, it is easily seen
 that
in an $H$-module algebra, grouplike elements act as unit
preserving automorphisms while primitive elements act as
derivations. In particular, for $H=kG$ the group algebra of a
discrete group, an $H$-module algebra structure on $A$ is simply
an action of $G$ by unit preserving automorphisms on $A$.
 Similarly, we
have a 1-1 correspondence between $U(\mathfrak{g})$-module algebra
structures on $A$ and Lie actions of the Lie algebra
$\mathfrak{g}$ by derivations on $A$.

  An algebra $B$ is called a left $H$-{\it comodule algebra}, if
$B$ is a left $H$-comodule and the multiplication and  unit maps of $B$ are $H$-comodule maps. That is
$$(ab)^{(-1)}\otimes (ab)^{(0)}= a^{(-1)} b^{(-1)}\otimes a^{(0)}b^{(0)}, \qquad (1_B)^{(-1)}\otimes
(1_B)^{(0)}=1_H\otimes 1_B,$$
for all $a, b$ in $B$.

  A left   $H$-{\it module coalgebra} is a coalgebra $C$ which is
 a left $H$-module such that   the comultiplication map $\Delta: C\to C\otimes C$ and the counit
 map $\varepsilon : C\to k$  are
 $H$-module maps. That is
$$(hc)^{(1)}\otimes (hc)^{(2)}=h^{(1)}c^{(1)}\otimes h^{(2)}c^{(2)}, \qquad \varepsilon (hc)=
\varepsilon (h)\varepsilon (c),$$
for all $h \in H$ and $c\in C$.

Finally,  a coalgebra $D$ is called a left $H$-{\it comodule coalgebra} when
the comultiplication and counit maps  of $D$ are morphisms of $H$-comodules. That is
$$ c^{(1)(-1)} c^{(2)(-1)}\otimes c^{(1)(0)}\otimes c^{(2)(0)}=c^{(-1)}\otimes
c^{(0)(1)}\otimes c^{(0)(2)},$$
$$\varepsilon (c)1_H=c^{(-1)}\varepsilon (c^{(0)}),$$
for all $c\in C.$

We call the above four types of symmetries, {\it symmetries of type A, B,  C,} and {\it D},
respectively. All four types of symmetries are present within an arbitrary Hopf
algebra.
For example, the coproduct $\Delta: H\to H\otimes H$ gives
 $H$ the structure of a left (and right) $H$-comodule algebra while the
product $H\otimes H\to H$ turns
  $H$ into  a left (and right) $H$-module coalgebra. These are noncommutative analogues of translation
  action of a group on itself. The
 {\it conjugation action} $H\otimes H \to H$, $g\otimes h \mapsto
 g^{(1)}hS(g^{(2)})$ gives $H$ the structure of a left $H$-module
 algebra. The {\it co-conjugation} action $H\to H\otimes H$, $ h \mapsto
 h^{(1)} S(h^{(3)})\otimes h^{(2)}$ turns $H$ into an $H$-comodule coalgebra.

 For a different example,
 we turn to the Connes-Moscovici Hopf algebra $\mathcal{H}_1$. An
 important feature of $\mathcal{H}_1$, and in fact its  {\it raison
 d'\^etre}, is that it acts as  quantum symmetries of  various objects of interest in noncommutative
 geometry, like the `space' of leaves of codimension one foliations or the `space' of modular forms
 modulo the action of Hecke correspondences. Let $M$ be a one
 dimensional manifold and $A=C^{\infty}_0(F^+M)$ denote the algebra of smooth
 functions with compact support on the bundle of positively oriented
 frames on $M$. Given a discrete group $\Gamma \subset Diff^+(M)$ of orientation preserving
 diffeomorphisms of $M$, one has a natural prolongation of the action of  $\Gamma$ to
  $F^+(M)$ by
  $$\varphi (y, y_1)=(\varphi (y), \varphi' (y)(y_1)).$$
Let $A_{\Gamma}= C^{\infty}_0(F^+M)\rtimes \Gamma$ denote the corresponding crossed product
algebra. Thus the elements of $A_{\Gamma}$ consist of finite linear combinations
(over $\mathbb{C}$)
of terms $fU_{\varphi}^{\ast}$ with $f\in C^{\infty}_0(F^+M)$ and  $\varphi \in \Gamma$. Its product
is defined by
$$fU_{\varphi}^{\ast} \cdot gU_{\psi}^{\ast}=(f \cdot \varphi (g))U_{\psi \varphi}^{\ast}.$$
 There is an
action of  $\mathcal{H}_1$ on $A_{\Gamma}$  given by \cite{cm2, cm5}:
$$Y(fU_{\varphi}^{\ast})= y_1 \frac{\partial f}{\partial y_1}
U_{\varphi}^{\ast}, \quad X(fU_{\varphi}^{\ast})= y_1 \frac{\partial f}{\partial y}
U_{\varphi}^{\ast},$$
$$\delta_n(fU_{\varphi}^{\ast})=y_1^n
\frac{ d^n} {dy^n} (log \frac{d\varphi}{dy})fU_{\varphi}^{\ast}.$$
Once these formulas are given, it can be checked, by a long
computation, that $A_{\Gamma}$ is indeed an $\mathcal{H}_1$-module
algebra. In the original application, $M$ is a transversal for a
codimension one foliation and thus $\mathcal{H}_1$ acts via
transverse differential operators \cite{cm2}.

We recall, very briefly,  the action of  the Hopf algebra
$\mathcal{H}_1$ on the so called {\it modular Hecke algebras},
discovered by Connes and Moscovici in \cite{cm5, cm6} where
detailed discussions and a very intriguing dictionary comparing transverse
geometry notions with modular forms notions can be found. For each
$N\geq 1$, let
$$\Gamma (N)= \left\{
\left( \begin{matrix} a & b\\c & d
\end{matrix} \right) \in SL(2, \mathbb{Z}); \; \; \left( \begin{matrix} a & b\\c & d
\end{matrix}\right)=\left( \begin{matrix} 1 & 0\\0 & 1
\end{matrix}\right) \;  \text{mod}\, N \right\}$$
denote the level $N$ congruence subgroup of $\Gamma (1)=SL(2,
\mathbb{Z})$. Let $\mathcal{M}_k (\Gamma (N))$ denote the space of modular
forms of level $N$ and weight $k$ and
$$\mathcal{M} (\Gamma (N)):=\oplus_k \mathcal{M}_k (\Gamma (N))$$
be the graded algebra of modular forms of level $N$. Finally let
$$\mathcal{M}:= \underset{\underset{N}{\to}}{\text{lim}} \; \mathcal{M} (\Gamma (N))$$
denote the algebra of modular forms of all levels, where the inductive system is defined by divisibility.
The group
$$G^{+}(\mathbb{Q}):= GL^+ (2, \mathbb{Q}),$$
acts  on  $\mathcal{M}$ by its usual action on functions on the upper half plane
(with corresponding weight):
$$(f, \alpha) \mapsto f|_k \alpha (z)= det(\alpha)^{k/2}(cz+d)^{-k}f(\alpha \cdot
z),$$
$$\alpha = \left( \begin{matrix} a & b\\c & d
\end{matrix}\right), \quad \alpha. z=\frac{az+b}{cz+d}.$$
The elements of the corresponding crossed-product algebra
$$\mathcal{A}=\mathcal{A}_{G^+(\mathbb{Q})}:= \mathcal{M}\rtimes G^+(\mathbb{Q}),$$
are finite sums
$$\sum fU^*_{\gamma}, \quad f\in \mathcal{M},\;\; \gamma \in
G^+(\mathbb{Q}),$$
with a product defined by
$$fU^*_{\alpha}\cdot gU^*_{\beta}=(f\cdot g|\alpha)U^*_{\beta \alpha}.$$
$\mathcal{A}$ can be thought of as the algebra of `noncommutative
coordinates' on the `noncommutative quotient space' of modular
forms modulo Hecke correspondences \cite{cm5}.

Consider the operator $X$ of degree two on the space of modular forms
defined by
$$X:= \frac{1}{2\pi i}\frac{d}{dz}-\frac{1}{12\pi i}\frac{d}{dz} (log
\Delta)\cdot Y,$$
where
$$\Delta (z)=(2\pi)^{12}\eta^{24}(z)=(2\pi)^{12}q\prod_{n=1}^{\infty} (1-q^n)^{24},
\quad q=e^{2\pi iz},$$
and $Y$ denotes the grading operator
$$Y(f)= \frac{k}{2}\cdot f, \quad \text{for all}\; f\in \mathcal{M}_k.$$
The following proposition is proved in \cite{cm5}. It shows that $\mathcal{A}_{G^+(\mathbb{Q}})$
is an $\mathcal{H}_1$-module algebra:
\begin{prop}
There is a unique action of $\mathcal{H}_1$ on
$\mathcal{A}_{G^+(\mathbb{Q})}$ determined by
$$ X(fU^*_{\gamma})=X(f)U^*_{\gamma}, \quad
Y(fU^*_{\gamma})=Y(f)U^*_{\gamma},$$
$$\delta_1( fU^*_{\gamma})=\mu_{\gamma}\cdot f(U^*_{\gamma}),$$
where
$$\mu_{\gamma} (z)=\frac{1}{2\pi
i}\frac{d}{dz} \text{log} \frac{\Delta|\gamma}{\Delta}.$$
\end{prop}

More generally, for any congruence subgroup $\Gamma$   an algebra
$A(\Gamma)$ is constructed in \cite{cm5} that contains  as
subalgebras both the algebra of $\Gamma$-modular forms and the
Hecke ring at level $\Gamma$. There is also a corresponding action
of $\mathcal{H}_1$ on $A(\Gamma)$.

For  symmetries of type  A, B, C, and D  there is a corresponding {\it crossed
product}, or {\it  smash product},
 construction that generalizes crossed products for actions of a
 group. We recall these constructions only for A and D symmetries,  as well as a more elaborate
 version called {\it bicrossed product} construction. We shall see that
 Connes-Moscovici's Hopf algebra $\mathcal{H}_1$ is a bicrossed product
 of two, easy to describe,
 Hopf algebras.

 Let $A$ be a left $H$-module algebra. The underlying vector space of the crossed product algebra
 $A\rtimes H$ is the vector space $A\otimes H$ and  its  product is determined by
 $$(a\otimes g)(b\otimes h)= a(g^{(1)}b)\otimes  g^{(2)}h.$$
 One checks that endowed with $1\otimes 1$ as its unit,  $A\rtimes H$ is an
 associative unital algebra. For example, for
 $H=kG$, the group algebra of a discrete group $G$ acting by automorphisms on an algebra
 $A$,  the algebra $A\rtimes \mathcal{H}$ is
  isomorphic to the crossed product algebra $A\rtimes G$.
For a second simple example, let a Lie algebra $\mathfrak{g}$ act
via derivations on a commutative algebra $A$. Then the crossed
product algebra $A \rtimes U(\mathfrak{g}) $ is a subalgebra of
the algebra of differential operators on $A$ generated by
derivations from $\mathfrak{g}$ and multiplication operators by
elements of $A$. The simples example is when $A=k[x]$ and
$\mathfrak{g}=k$ acting via $\frac{d}{dx}$ on $A$. Then $A\rtimes
U(\mathfrak{g})$ is the {\it weyl algebra} of differential
operators on the line with polynomial coefficients.

Let $D$ be a right $H$-comodule coalgebra via the coaction $
d\mapsto d^{({0})}\ot d^{({1})}\in D\ot H$. The underlying
linear space of the {\it crossed product coalgebra} $H\rtimes D$
is $H\otimes D$. It is a coalgebra whose  coproduct and counit are defined by
$$ \Delta (h\otimes d)= h^{(1)}\ot (d^{(1)})^{({0})}\ot h^{(2)}(d^{(1)})^{({1})}\ot d^{(2)},
\quad \varepsilon (h\otimes d)=  \varepsilon
(d) \varepsilon (h).$$

The above two constructions deform multiplication or comultiplication, of algebras  or coalgebras,
respectively. Thus to obtain
a simultaneous deformation of multiplication and comultiplication of a
Hopf algebra it stands to reason to try to apply both constructions simultaneously. This idea, going back
to G. I. Kac in 1960's in the context of Kac-von Neumann Hopf algebras,
has now found its complete generalization in the notion of {\it bicrossed
product}  of {\it matched pairs} of Hopf algebras. See Majid's book \cite{maj}
for extensive  discussions and references. There are many variations of
this construction of which the most relevant for the structure of
the Connes-Moscovici Hopf algebra is the following.

Let $U$ and $F$ be two Hopf algebras. We assume that $F$ is a
 left $U$-module algebra and $U$ is a right  $F$-comodule coalgebra via $\rho:U\rightarrow U\ot F$. We say
that $ (U, F)$ is a {\it matched pair} if the action and coaction
satisfy the compatibility condition:
\begin{align*}
\epsilon(u(f))=\epsilon(u)\epsilon(f),&
~~~\Delta(u(f))=(u^{(1)})^{(0)}(f^{(1)})\ot
(u^{(1)})^{({1})}(u^{(2)}(f^{(2)})),\\
\rho(1)=1\ot 1,& ~~\rho(uv)=(u^{(1)})^{({0})} v^{({0})}\ot
(u^{(1)})^{({1})}(u^{(2)}(v^{({1})})),\\
(u^{(2)})^{({0})}\ot &(u^{(1)}(f))(u^{(2)})^{({1})}
=(u^{(1)})^{({0})}\ot (u^{(1)})^{({1})}(u^{(2)}(f)).
\end{align*}
 Given a matched pair as above, we define its bicrossed product
Hopf algebra $F\rtimes U$ to be $F\otimes U$ with crossed product
algebra structure and crossed coproduct coalgebra structure. Its
antipode $S$ is defined by
$$S(f\ot u)=(1\ot S(u^{({0})}))(S(fu^{({1})})\ot 1).$$
It is a remarkable fact that, thanks to the the above compatibility conditions, all the axioms of a Hopf
algebra are satisfied for $F\rtimes U$.

The simplest and first example of this bicrossed product construction
is as follows. Let $G=G_1 G_2$ be a {\it factorization} of a finite group
$G$. This means that $G_1$ and $G_2$ are subgroups of $G$,
$G_1\cap G_2=\{e\}, $ and $G_1 G_2=G$. We denote the factorization of $g$ by $g=g_1g_2$. The relation
$g\cdot h:=(gh)_2$ for $ g\in G_1$ and $h\in G_2$  defines a left action
of $G_1$ on $G_2$. Similarly $g\bullet h:=(gh)_1$ defines a right action
of $G_2$ on $G_1$. The first action turns $F=F(G_2)$ into a left $U=kG_1$-module algebra. The second
action turns $U$ into a right $F$-comodule coalgebra. The later coaction is simply the dual
of the map $F(G_1) \otimes kG_2 \to F(G_1)$ induced by the right action
of
$G_2$ on $G_1$. Details of this example can be found in \cite{maj} and \cite{cm2}.

\begin{example}
1. By a theorem of  Kostant \cite{sw},  any cocommutative Hopf algebra $H$ over an
algebraically
closed field $k$ of characteristic zero is isomorphic (as  a Hopf algebra) with a
  crossed product algebra  $H= U(P(H))\rtimes  k G(H)$,
where  $P(H)$ is the  Lie algebra  of  primitive elements of $H$ and $G(H)$ is the
 group of all grouplike
 elements of $H$ and $G(H)$ acts on $P(H)$ by inner
 automorphisms $(g,h)\mapsto ghg^{-1}$, for $g\in G(H)$ and $h\in P(H)$.
 The coalgebra structure of $H= U(P(H))\rtimes  k G(H)$ is simply the
 tensor product of the two coalgebras  $U(P(H))$ and $  k G(H).$\\

\noindent 2. We show that the Connes-Moscovici Hopf algebra is a
bicrossed Hopf algebra. Let $G=Diff (\mathbb{R})$ denote the group of
diffeomorphisms of the real line. It has a factorization of the form
$$G=G_1G_2,$$
where $G_1$ is the subgroup of diffeomorphisms that satisfy
 $$\varphi (0)=0, \quad \varphi'(0)=1,$$
and $G_2$ is the $ax+b$- group of affine transformations.  The first
Hopf algebra, $F$, is
formally speaking, the {\it algebra of polynomial functions} on the prounipotent group
$G_1$. It can also be defined  as the ``continuous dual'' of the
enveloping algebra of the Lie algebra of $G_1$. It is a commutative Hopf
algebra generated by functions $\delta_n$, $n=1,2,\dots$, defined by
$$\delta_n (\varphi)=\frac{d^n}{dt^n}(log (\varphi'(t))|_{t=0}.$$
  The second Hopf algebra,
$U$,  is the universal enveloping algebra of the Lie algebra
$\mathfrak{g}_2$ of the $ax+b$-group. It has generators $X$ and $Y$ and one relation $[X,
Y]=X$.

 The Hopf algebra $F$ has a right  $U$-module algebra structure
 defined by
 $$\delta_n (X)=-\delta_{n+1}, \quad \text{ and}\;\;
\delta_n (Y)=-n\delta_n.$$
 The Hopf algebra $U$ is a left
$F$-comodule coalgebra via
$$X\mapsto 1\ot X+ \delta_1\ot X, \quad \text{ and}\;\;
Y\mapsto 1\ot Y.$$ One can check that they are a matched pair of  Hopf
algebras and the resulting bicrossed product Hopf algebra is
isomorphic to the Connes-Moscovici Hopf algebra $\mathcal{H}_1$. See \cite{cm2} for a slightly different
approach and fine points of the proof. \\

\noindent 3.  Another important example of a bicrossed construction is the
Drinfeld double $D(H)$ of a finite dimensional Hopf algebra $H$ defined
as a bicrossed product $H\rtimes H^*$ \cite{maj}.
\end{example}

\subsection{Modular pairs}
Let $H$ be a Hopf algebra, $\delta: H\to k$ a character and
$\sigma \in H$ a grouplike element. Following \cite{cm2, cm3}, we
say that $(\delta, \sigma)$ is a {\it modular pair} if $\delta
(\sigma)=1$, and a {\it modular pair in involution} if in addition
we have
$$\widetilde{S}_{\delta}^2=Ad_{\sigma}, \quad \text{or},  \quad \widetilde{S}_{\delta}^2 (h)=\sigma
h\sigma^{-1},$$
for all $h$ in $H$. Here the
  $\delta$-{\it twisted antipode}
$\widetilde{S}_{\delta}: H \rightarrow H$ is defined by $\widetilde{S}_{\delta}=\delta *
S$, i.e.
$$\widetilde{S}_{\delta}(h)= \delta (h^{(1)})S(h^{(2)}),$$
for all $h \in H$.

The notion of an invariant trace for actions of groups and Lie
algebras can be extended to the Hopf setting. For applications to
transverse geometry and number theory, it is important to
formulate a notion of `invariance' under the presence of a modular
pair.  Let $A$ be an $H$-module algebra, $\delta $ a character of
$H$, and $\sigma \in H$ a grouplike element.  A $k$-linear map
$\tau :A \rightarrow k$ is called $\delta$-{\it invariant} if for
all $h \in H$ and $a\in A$,
$$ \tau (h(a))=\delta (h) \tau (a).$$
$\tau$ is called a $\sigma$-{\it trace} if for all $a, b$ in $A$,
$$ \tau (ab)=\tau (b \sigma (a)).$$
For the following lemma from \cite{cm4} the fact
that $A$ is unital is crucial. For $a, b \in A,$ let
$$<a, \, b>:=\tau (ab).$$
\begin{lem} (Integration by parts formula). Let $\tau $ be a $\sigma$-trace on $A$. Then $\tau$ is
$\delta$-invariant if and only if
the {\it integration by parts formula} holds:
\begin{eqnarray*}
<h(a), \, b>=<a, \,  \widetilde{S}_{\delta}(h)(b)>,
\end{eqnarray*}
for all $h \in H$ and $a, b \in A$.
\end{lem}

Loosely speaking, the lemma says that the formal adjoint of the differential operator $h$ is
$\widetilde{S}_{\delta}(h)$.

\begin{example} 1. For any Hopf algebra $H$, the pair $(\varepsilon, 1)$ is modular. It is involutive if and
only if $S^2=id$. This happens, for example, when $H$ is a commutative or cocommutative Hopf algebra.\\

\noindent 2. The original non-trivial example of a modular pair in involution
is the pair $(\delta, 1)$ for
Connes-Moscovici  Hopf algebra $\mathcal{H}_1$. Let $\delta $
denote the unique extension of the modular character
$$\delta : \mathfrak{g}_{aff} \to \mathbb{R}, \quad \delta(X)=1, \; \delta(Y)=0,$$
to a character $\delta : U (\mathfrak{g}_{aff})\to \mathbb{C}.$
There is a unique extension of $\delta $ to a character, denoted
by the same symbol  $\delta : \mathcal{H}_1 \to \mathbb{C}.$
Indeed the relations $[Y, \delta_n]=n\delta_n$ show that we must
have $\delta (\delta_n)=0$, for $n=1, 2, \cdots.$ One can then
check that these relations are compatible with the algebra
structure of $\mathcal{H}_1.$ Recall the algebra
$A_{\Gamma}=C^{\infty}_0(F^+(M)\rtimes \Gamma$ from Section 1.2.
It admits a $\delta$-invariant trace $\tau : A_{\Gamma} \to
\mathbb{C}$ under its canonical $\mathcal{H}_1$ action given by
\cite{cm2}:
$$\tau (fU^{\ast}_{\varphi})=\int_{F^+(M)}f(y, y_1)\frac{dy dy_1}{y_1^2},
\quad  \text{if} \; \varphi =1,$$
and $\tau (fU^{\ast}_{\varphi})=0$, otherwise.

\noindent 3. Let $H=A(SL_q(2,k))$ denote the Hopf algebra of functions on quantum $SL(2, k)$.
 As an algebra it is generated by symbols $x,\; u,\; v,\; y,$ with the following relations:
$$ ux=qxu, \;\; vx=qxv, \;\; yu=quy,\;\;yv=qvy,$$
$$uv=vu,\;\; xy-q^{-1}uv=yx-quv=1. $$
The coproduct, counit and antipode of   $\mathcal{H}$ are defined by
$$\Delta (x)=x \otimes x+u \otimes v,\;\;\;\Delta (u)=x \otimes u+u \otimes y, $$
$$\Delta (v)=v \otimes x+y \otimes v,\;\;\;\Delta (y)=v \otimes u+y \otimes y, $$
$$\epsilon (x)=\epsilon (y)=1,\;\;\;\epsilon (u)=\epsilon (v)=0,$$
$$S(x)=y,\;\; S(y)=x,\;\;S(u)=-qu,\;\;S(v)=-q^{-1}v. $$
Define a character  $\delta : H \to k $ by:
$$\delta (x)=q,\; \; \delta(u)=0,\; \; \delta (v)=0,\; \; \delta (y)=q^{-1}.$$
One checks that  $\widetilde {S}_{\delta}^2=id$. This shows that $(\delta, 1)$ is a modular pair for
$H$. This example and its Hopf-cyclic cohomology is studied in
\cite{kr1}.

More generally, it is shown in \cite{cm3} that {\it coribbon Hopf algebras} and compact
quantum groups are endowed with canonical modular pairs of the form $(\delta, 1)$
and,
dually,  ribbon Hopf algebras have canonical modular pairs of the type
$(1, \sigma)$.
\\

\noindent 4. We will see in the next section that modular pairs in
involution are in fact
 one dimensional cases of what we call stable anti-Yetter-Drinfeld modules, i.e. they are
 one dimensional
noncommutative local systems that one can introduce into Hopf algebra equivariant cyclic cohomology.

\end{example}

\section{Anti-Yetter-Drinfeld modules}

An important question left open in our paper \cite{kr2} was the
    issue of identifying the most general class of coefficients allowable in cyclic (co)homology of
    Hopf
    algebras and Hopf-cyclic cohomology in general. This
    problem
    is  completely solved, among other things, in \cite{hkrs1}. It is shown in this paper that
    the most general coefficients are the
    class of  so called  stable anti-Yetter-Drinfeld modules. In Section 3.2 we briefly
    report on this very recent development as well.

It is quite surprising that  when the general formalism of cyclic
cohomology theory, namely the theory of (co)cyclic modules
\cite{c1},  is applied to Hopf algebras, variations of  such
standard notions like Yetter-Drinfeld (YD) modules appear
naturally. The so called anti-Yetter-Drinfeld modules introduced
in \cite{hkrs1} are twistings, by modular pairs in
involution, of YD modules.  This means that the category of
anti-Yetter-Drinfeld modules is a ``mirror image'' of the category
of YD modules. We mention that anti-Yetter-Drinfeld modules were
considered independently also by C. Voigt in connection with his work on
equivariant cyclic cohomology \cite{voi}.

\subsection{Yetter-Drinfeld modules}
Yetter-Drinfeld modules were introduced by D. Yetter under the name
{\it crossed bimodules} \cite{y}. The present name was coined in
\cite{rt}. One of the  motivations was to define a {\it braiding} on the
monoidal category  $H-Mod$ of representations of a, not necessarily commutative or
cocommutative, Hopf algebra $H$.
To define  such a braiding one should either restrict to special  classes of
Hopf algebras, or, to special classes of modules.   Drinfeld showed that
when $H$ is a
{\it quasitriangular} Hopf algebra, then $H-Mod$ is a
braided monoidal category. See \cite{maj} for definitions and references.
Dually, when $H$ is {\it coquasitriangular},
the category $H-Comod$ of $H$-comodules is a braided monoidal category.
In \cite{y} Yetter shows that to obtain a braiding on a subcategory of $H-Mod$, for an
arbitrary $H$, one has essentially one choice and that is restricting to
the class of Yetter-Drinfeld modules as we explain now.

Let $H$ be a Hopf algebra and $M$ be a left $H$-module and left $H$-comodule simultaneously. We say
that $M$ is a left-left {\it Yetter-Drinfeld $H$-module} if  the two structures
on $M$ are compatible in the sense that
$$\rho(hm)=h^{(1)}m^{(-1)}S(h^{(3)})\ot h^{(2)}m^{(0)},$$
 for all $h\in H$ and $m\in M$ \cite{maj, rt, y}. One can similarly define left-right, right-left and
 right-right YD modules. A morphism of YD modules $M \to N$ is an $H$-linear and $H$-colinear map
 $f:M\to N$.  We denote the category of left-left YD modules over $H$ by
 $^{H}_{H} \mathcal{Y}\mathcal{D}$.

 This notion is closely related to the
 Drinfeld double of finite dimensional Hopf algebras. In fact if $H$ is
 finite dimensional, then one  can show that the category $^{H}_{H} \mathcal{Y}\mathcal{D}$  is
 isomorphic to the category of left modules over the Drinfeld double
 $D(H)$ of $H$ \cite{maj}.

The following facts about the category $^{H}_{H} \mathcal{Y}\mathcal{D}$ are well known \cite{maj, rt, y}:

1.  The tensor product $M\otimes N$ of two YD modules is a YD module. Its module and comodule structure
are the standard ones:
$$h (m\otimes n)=h^{(1)}m\otimes h^{(1)}n, \quad (m\otimes n)\mapsto
m^{(-1)}n^{(-1)}\otimes m^{(0)}\otimes n^{(0)}.$$
 This shows that the category  $^{H}_{H} \mathcal{Y}\mathcal{D}$ is a monoidal subcategory  of the monoidal
 category $H-Mod$.\\

2. The category $^{H}_{H} \mathcal{Y}\mathcal{D}$ is a braided monoidal category under the
braiding
$$c_{M, N}: M\otimes N \to N\otimes M, \quad m\otimes n \mapsto m^{(-1)}\cdot n\otimes  m^{(0)}.$$
In fact Yetter proves a strong inverse to this statement as well \cite{y}: for any
small strict monoidal category $\mathcal{C}$ endowed with  a monoidal
functor $F: \mathcal{C}\to Vect_f$ to the category of finite dimensional
vector spaces, there is a Hopf algebra $H$ and a monoidal functor
$\tilde{F}: \mathcal{C} \to  ^{H}_{H} \mathcal{Y}\mathcal{D}$ such that the following diagram
comutes\\

\begin{center}
$\begin{xy}
\xymatrix{ \mathcal{C}\ar[r]\ar[d]&  ^{H}_{H} \mathcal{Y}\mathcal{D} \ar[d]\\
Vect_f\ar[r]& Vect}
\end{xy}$
\end{center}

3. The category $^{H}_{H} \mathcal{Y}\mathcal{D}$ is the {\it
center} of the monoidal category $H-Mod$. Recall that the (left)
center $\mathcal{Z}\mathcal{C}$ of a monoidal category \cite{kas} is
a category whose objects are pairs $(X, \sigma_{X, -})$, where $X$
is an object of $\mathcal{C}$ and $\sigma_{X, -}: X\otimes -\to
-\otimes X$ is a natural isomorphism satisfying certain
compatibility axioms with associativity and unit constraints. It
can be shown that the center of a monoidal category is a braided
monoidal category and
$\mathcal{Z} (H-Mod)=^{H}_{H} \mathcal{Y}\mathcal{D}$  \cite{kas}. \\

\begin{example}
1. Let $H=kG$ be the group algebra of a discrete group $G$.  A left $kG$-comodule is
simply a $G$-graded vector space
$$M=\bigoplus_{g\in G}M_g$$
 where
the coaction is defined by
$$m\in M_g \mapsto g\otimes m.$$
 An action of
$G$ on $M$ defines a YD module structure iff for all $g, h \in G$,
$$hm\in M_{hgh^{-1}}.$$
This example can be  explained as
follows. Let $\mathcal{G}$ be a groupoid whose objects are $G$ and its morphisms are defined by
$$Hom(g, h)=\{k\in G; \; kgk^{-1}=h\}.$$
 Recall that an {\it action}
 of a groupoid $\mathcal{G}$ on the category $Vect$ of vector spaces is simply a
 functor $F: \mathcal{G} \to Vect$.
  Then it is easily seen that we have a
one-one correspondence between YD modules for $kG$ and actions of $\mathcal{G}$ on $Vect$. This example
clearly shows that one can think of an YD module over $kG$ as an {\it `equivariant sheaf'} on $G$ and
of $YD$ modules as noncommutative analogues of equivariant sheaves on a topological  group.\\

\noindent 2. If $H$ is cocommutative then any left $H$-module $M$
can be turned into a left-left YD module  via the coaction
$m\mapsto 1\otimes m$. Similarly, when $H$ is cocmmuntative then
any left $H$-comodule $M$ can be turned into a YD module
via the $H$-action $h\cdot m:=\varepsilon(h)m.$\\

\noindent 3. Any Hopf algebra acts on itself via {\it conjugation action} $g \cdot h:=g^{(1)}hS(g^{(2)})$ and coacts
via translation coaction $h\mapsto h^{(1)}\otimes h^{(2)}.$ It can be
checked that this endows $M=H$ with a YD module structure.
\end{example}

\subsection{Stable anti-Yetter-Drinfeld modules}
This class of modules for Hopf algebras were introduced for the
first time in \cite{hkrs1}. Unlike Yetter-Drinfeld modules, its
definition, however, was entirely motivated and dictated by cyclic
cohomology theory: the  anti-Yetter-Drinfeld condition guarantees
that the simplicial and cyclic operators are well defined on
invariant  complexes and the stability condition implies
that the crucial periodicity condition for cyclic modules are
satisfied.

\begin{definition} A left-left  anti-Yetter-Drinfeld  $H$-module is a left $H$-module and
left $H$-comodule such that
$$\rho(hm)=h^{(1)}m^{(-1)}S(h^{(3)})\ot h^{(2)}m^{(0)},$$
for all $h\in H$ and  $m\in M.$
We say that $M$ is stable if in addition we have
$$m^{(-1)}m^{(0)}=m,  $$
for all $m\in M$.
\end{definition}

There are of course analogous definitions for left-right,
right-left and right-right stable anti-Yetter-Drinfeld (SAYD)
modules. We note that by changing $S$ to $S^{-1}$ in the above
equation, we obtain the compatibility condition for a
Yetter-Drinfeld  module from the previous subsection.

The following lemma from \cite{hkrs1} shows that 1-dimensional
SAYD modules correspond to Connes-Moscovici's modular pairs in
involution:
\begin{lem} There is a one-one correspondence between modular pairs in
involution $(\delta, \sigma)$ on $H$ and SAYD module structure on $M=k$, defined by
$$h. r=\delta (h) r, \quad r\mapsto \sigma \otimes r,$$
for all $h\in H$ and $r\in k$. We denote this module by $M=^{\sigma}\!\!\!k_{\delta}.$
\end{lem}

Let $^{H}_{H} \mathcal{A}\mathcal{Y}\mathcal{D}$ denote the
category of left-left anti-Yetter-Drinfeld $H$-modules, where
morphisms are $H$-linear and $H$-colinear maps. Unlike YD modules,
anti-Yetter-Drinfeld modules do not form a monoidal category under
the standard tensor product. This can be checked easily on
1-dimensional modules given by modular pairs in involution.
 The following result of  \cite{hkrs1}, however,
shows that the tensor
product of an anti-Yetter-Drinfeld module with a Yetter-Drinfeld module
is again anti-Yetter-Drinfeld.
\begin{lem}
Let $M$ be a Yetter-Drinfeld module and   $N$ be an anti-Yetter-Drinfeld module (both left-left).
  Then $M \otimes N $ is an anti-Yetter-Drinfeld module under the
  diagonal action and coaction:
  $$h (m\otimes n)=h^{(1)}m\otimes h^{(1)}n, \quad (m\otimes n)\mapsto
m^{(-1)}n^{(-1)}\otimes m^{(0)}\otimes n^{(0)}.$$
\end{lem}

In particular, using a modular pair in involution $(\delta, \sigma)$, we
obtain a functor
$$^{H}_{H} \mathcal{Y}\mathcal{D} \to ^{H}_{H} \mathcal{A}\mathcal{Y}\mathcal{D},
\quad M\mapsto \overset{-}{M}= ^{\sigma}\!\!\!k_{\delta}\otimes M.$$

This result clearly shows that AYD modules can be regarded as the
{\it twisted analogue}  or {\it mirror image} of YD modules, with
twistings provided by modular pairs in involution. This result was
later strengthened by the following result,  pointed out to us by
M. Staic \cite{sta}. It shows that if the Hopf algebra has a
modular pair in involution then the  category of YD modules is equivalent to the category of AYD
modules:
 \begin{prop}
 Let H be a Hopf algebra, $(\delta, \sigma)$  a modular pair in involution
and M an anti-Yetter-Drinfeld module. If we define
$m \cdot h =
mh^{(1)} \delta (S(h^{(2)}))$  and $\rho(m) = \sigma^{-1} m^{(-1)}\otimes m^{(0)}$,  then
$(M, \cdot, \rho)$ is an Yetter-Drinfeld module. Moreover in this way we get an isomorphism
between the categories of AYD modules and YD mpdules.
\end{prop}
It follows that  tensoring
with $^{\sigma^{-1}}\!\!\!k_{\delta \circ S}$
 turns the anti-Yetter-Drinfeld modules to Yetter-Drinfeld modules and this   is the inverse for
 the operation of tensoring with $^{\sigma}\!\!\!k_{\delta}$.

\begin{example}
1. For Hopf algebras with $S^2=I$, e.g. commutative or cocommutative Hopf algebras, there
is no distinction between YD and AYD modules. This applies in particular
to $H=kG$ and Example 2.1.1. The stability condition
$m^{(-1)}m^{(0)}=m$ is equivalent to
$$g\cdot m=m, \quad \text{for all} \, \, \, g\in G, \, m\in M_g.$$\\

\noindent 2. {\it Hopf-Galois extensions} are noncommutative
analogues of principal bundles in (affine) algebraic geometry.
Following \cite{hkrs1} we  show that they give rise to large
classes of examples of SAYD modules. Let $P$ be a right
$H$-comodule algebra, and let
$$B:=P^H=\{p\in P; \, \rho(p)=p\otimes 1\}$$
be the space of coinvariants of $P$. It is easy to see that $B$ is
a subalgebra of $P$. The extension $B\subset P$ is called a
Hopf-Galois extension if the map
$$can:  P\otimes_B P\to B\otimes H, \quad p\otimes p'\mapsto p\rho (p'), $$
is bijective. (Note that in the commutative case this corresponds
to the condition that the action of the structure group on fibres
is free). The bijectivity assumption allows us to
define the translation map $T: H \rightarrow P\ot_B P$,
$$T(h)=can^{-1}(1\ot h)=h^{(\bar{1})}\ot h^{(\bar{2})}.$$
It can be
checked that
 the centralizer  $Z_B(P)=\{p\mid bp=pb~~ \forall b\in B\}$ of $B$
 in $P$ is a subcomodule of $P$. There is an action of $H$ on
 $Z_B(P)$ defined by $ph=h^{({1})} p h^{({2})}$ called the
 Miyashita-Ulbrich action. It is shown that this action and
 coaction satisfy the Yetter-Drinfeld compatibility condition. On
 the other hand if $B$ is central, then by defining the new action $ph=(S^{-1}(h))^{({2})} p
 (S^{-1}(h))^{({1})}$and the right coaction of $P$ we have a SAYD
 module. This example was the starting point of \cite{js} where the relative cyclic homology
 of Hopf-Galois extensions is related to a variant of Hopf-cyclic cohomology with
 coefficients in stable
 anti-Yetter-Drinfeld modules.
 In \cite{kr4} we showed that the  theory introduced in \cite{js} is ismorphic to one of the
 theories introduced in \cite{hkrs2}. \\

\noindent 3. Let $M=H$. Then with conjugation action $g\cdot h=
g^{(1)}hS(g^{(2)})$ and comultiplication $h\mapsto h^{(1)}\otimes
h^{(2)}$ as coaction, $M$ is an SAYD module.
\end{example}

 \section{Hopf-cyclic cohomology}

In this section we first recall the approach by  Connes and Moscovici
towards the definition of their
cyclic cohomology theory for Hopf algebras. The characteristic map
$\chi_{\tau}$ palys an imporant role here. Then we switch to the point
of view adopted in \cite{kr2} based on invariant complexes, culminating
in our joint work \cite{hkrs1, hkrs2}. The resulting Hopf-cyclic
cohomology theories of type $A$, $B$, and $C$, and their corresponding cyclic modules,
 contain all known examples of cyclic theory discovered so far. We note
 that very recently  A. Kaygun has extended the Hopf-cyclic cohomology to
 a cohomology for bialgebras with coefficients in stable modules. For Hopf algebras it reduces
 to Hopf-cyclic cohomology \cite{kay}. 

\subsection{Connes-Moscovici's  breakthrough}

Without going into  details, in this section we  formulate one of
 the problems that was faced and solved by Connes and Moscovici in the
course of their study of an index problem on foliated manifolds
\cite{cm2}. See also \cite{cm4} for a survey. As we shall see, this led them to a new cohomology
theory for Hopf algebras that is now an important example of  Hopf-cyclic
cohomology.

The local
index formula of Connes and Moscovici \cite{cm1} gives the Chern character
$Ch (A, h, D)$ of a regular
spectral triple $(A, h, D)$ as a  cyclic cocycle in the $(b, B)$-bicomplex of the algebra $A$.
For  spectral
triples of interest in transverse geometry \cite{cm2}, this cocycle is {\it differentiable} in
the sense that it is in the
image of the Connes-Moscovici characteristic map $\chi_{\tau}$ defined below, with
$H=\mathcal{H}_1$ and $A=\mathcal{A}_{\Gamma}$. To identify
this class in terms of characteristic classes of foliations, it would be
extremely helpful to
 show that it is the image of a cocycle for a cohomology theory for Hopf algebras. This is rather
similar to the situation for classical characteristic classes which are pull backs of group
cohomology classes.

We can formulate this problem abstractly as follows: Let $H$ be
a Hopf algebra endowed with a modular pair in involution $(\delta, \sigma)$,
and $A$ be an $H$-module algebra. Let $\tau: A \to
k$ be a $\delta$- invariant $\sigma$-trace on $A$ as defined in Section 1.3.
 Consider the Connes-Moscovici {\it
characteristic map}
$$\chi_{\tau}: H^{\otimes n}\longrightarrow Hom (A^{\otimes (n+1)},\;
k),$$
\begin{eqnarray*}
\chi_{\tau}(h_1\otimes \cdots \otimes h_n)(a_0\otimes \cdots \otimes a_n)=
\tau (a_0 h_1(a_1)\cdots h_n(a_n)).
\end{eqnarray*}
Now the burning question  is: can we   promote the collection of
spaces $\{H^{\otimes n}\}_{n\geq 0}$ to a {\it cocyclic module}
such that the characteristic map $\chi_{\tau}$ turns into a
morphism of cocyclic modules? We recall that the face, degeneracy,
and cyclic operators for  $Hom (A^{\otimes (n+1)},\; k) $ are
defined by \cite{c1}:
\begin{eqnarray*}
\delta^n_i \varphi (a_0, \cdots, a_{n+1})&=&\varphi (a_0, \cdots,
a_ia_{i+1}, \cdots, a_{n+1}), \quad i=0, \cdots, n,\\
\delta^n_{n+1} \varphi (a_0, \cdots, a_{n+1})&=&\varphi (a_{n+1}a_0,
a_1, \cdots, a_n),\\
\sigma^n_i\varphi (a_0, \cdots, a_{n}) &=&\varphi (a_0, \cdots,a_i, 1,
\cdots, a_{n}),\quad i=0, \cdots, n,\\
\tau_n \varphi (a_0, \cdots, a_n)& =&\varphi (a_n, a_0, \cdots, a_{n-1}).
\end{eqnarray*}

 The
relation
$$h(ab)=h^{(1)}(a)h^{(2)}(b)$$
shows that  in order for  $\chi_{\tau}$ to be compatible with  face
operators, the face operators
on $H^{\otimes n}$ must involve the coproduct of
$H$. In fact if we define, for $ 0\leq i \leq n$,
$\delta^n_i: H^{\otimes n} \to H^{\otimes (n+1)}$, by
\begin{eqnarray*}
\delta^n_0(h_1\otimes
\cdots \otimes h_n)&=& 1\otimes h_1\otimes
\cdots \otimes h_n,\\
 \delta^n_i(h_1\otimes
\cdots \otimes  h_n)&=& h_1\otimes \cdots \otimes h_i^{(1)}\otimes h_i^{(2)}\otimes \cdots
\otimes h_n,\\
\delta^n_{n+1}(h_1\otimes \cdots \otimes h_n)&=&h_1\otimes
\cdots \otimes  h_n\otimes \sigma,
\end{eqnarray*}
then we have, for all $n$ and $i$,
$$\chi_{\tau}\delta_i^n =\delta_i^n \chi_{\tau}.$$

 Similarly, the relation $h(1_A)=\varepsilon (h) 1_A,$
shows that the degeneracy operators on $H^{\otimes n}$ should
involve the counit of $H$. We thus define
$$\sigma^n_i(h_1\otimes \cdots \otimes h_n)=h_1 \otimes \cdots \otimes \varepsilon
(h_i)\otimes \cdots \otimes h_n.$$

The most difficult part in this regard is to guess the form of the
 {\it cyclic operator} $\tau_n : H^{\otimes n} \to H^{\otimes n}$.
 Compatibility with $\chi_{\tau}$ demands that
$$ \tau (a_0  \tau_n(h_1\otimes \cdots \otimes h_n)(a_1\otimes \cdots \otimes a_n))=
\tau (a_n  h_1(a_0)h_2(a_1)\cdots h_n(a_{n-1})),$$
for all $a_i$'s and $h_i$'s.
Now  integration by parts formula in Lemma 2.1, combined with the $\sigma$-trace property of $\tau$,
gives us:
$$\tau (a_1 h(a_0))=\tau (h(a_0)\sigma (a_1))=\tau (a_0\tilde{S}_{\delta}(h)(\sigma (a_1)).$$
This suggests that we should define $\tau_1: H \to H$ by
$$\tau_1(h)= \tilde{S}_{\delta}(h)\sigma.$$
Note that the condition
$\tau_1^2=I$ is equivalent to the involutive condition
$\tilde{S}_{\delta}^2=Ad_{\sigma}$.

For any $n$,
integration by parts formula together with the $\sigma$-trace property shows
that:
\begin{eqnarray*}
\tau (a_n h_1(a_0) \cdots h_n(a_{n-1}))& =&\tau (
h_1(a_0) \cdots h_n(a_{n-1})\sigma (a_n))\\
&=&\tau(a_0 \tilde{S}_{\delta}(h_1)(h_2 (a_1)\cdots h_n(a_{n-1})\sigma (a_n)))\\
&=& \tau(a_0 \tilde{S}_{\delta}(h_1)\cdot(h_2 \otimes \cdots \otimes h_n\otimes \sigma)
(a_1\otimes \cdots \otimes a_n).
\end{eqnarray*}
This suggests that the {\it Hopf-cyclic operator} $\tau_n : H^{\otimes
n}\to  H^{\otimes n}$ should be defined as
$$\tau_n (h_1\otimes \cdots \otimes h_n)=\tilde{S}_{\delta}(h_1)\cdot(h_2 \otimes \cdots \otimes
h_n\otimes \sigma),$$
where $\cdot$ denotes the diagonal action defined by
$$h\cdot (h_1\otimes \cdots \otimes h_n):= h^{(1)}h_1\otimes
h^{(2)}h_2\otimes \cdots \otimes h^{(n)}h_n.$$ We let $\tau_0=I:
H^{\otimes 0}=k \to H^{\otimes 0}$, be the identity map. The
remarkable fact, proved by Connes and Moscovici \cite{cm2, cm3},
is that endowed with the above face, degeneracy, and cyclic
operators, $\{H^{\otimes n} \}_{n\geq 0}$ is a cocyclic module.
The resulting cyclic cohomology groups are denoted by
$HC^n_{(\delta, \sigma)}(H)$, $n=0, 1,\cdots$ and we obtain the
desired  characteristic map
$$\chi_{\tau}: HC^n_{(\delta, \sigma)}(H) \to HC^n(A).$$

As with any cocyclic module, cyclic cohomology  can also be
defined in terms of cyclic cocycles. In this case a cyclic n-cocycle is
an element $x\in H^{\otimes n}$ satisfying the conditions
$$bx=0, \quad (1-\lambda)x=0,$$
where $b: H^{\otimes n} \to H^{\otimes (n+1)}$ and $\lambda:  H^{\otimes
n}\to H^{\otimes n}$ are defined by
\begin{eqnarray*}
 b(h^1\otimes \cdots \otimes h^n)&=&1\otimes h_1\otimes
\cdots \otimes h_n\\
&+&\sum_{i=1}^n (-1)^i  h_1\otimes \cdots \otimes h_i^{(1)}\otimes h_i^{(2)}\otimes \cdots
\otimes h_n\\
&+&(-1)^{n+1}h_1\otimes
\cdots \otimes  h_n\otimes \sigma, \\
\lambda (h_1\otimes \cdots \otimes h_n)&=& (-1)^n\tilde{S}_{\delta}(h_1)\cdot(h_2 \otimes \cdots \otimes
h_n\otimes \sigma).
\end{eqnarray*}

The cyclic cohomology groups $HC^n_{(\delta, \sigma)}(H)$ are
computed in several cases in \cite{cm2}. Of particular interest
for applications to transverse index theory and number theory is
the the (periodic)  cyclic cohomology of the Connes-Moscovici Hopf
algebra $\mathcal{H}_1$. It is shown in \cite{cm2} that the
periodic groups $HP^n_{(\delta, 1)}(\mathcal{H}_1)$ are canonically
isomorphic to the Gelfand-Fuchs cohomology of the Lie algebra of
formal vector fields on the line:
$$H^*(\mathfrak{a}_1, \mathbb{C})=HP^*_{(\delta,
1)}(\mathcal{H}_1).$$
Calculation of
the unstable groups is an interesting open problem. The following
interesting elements are however already been identified.
It can be directly checked that the elements $\delta_2' :
=\delta_2-\frac{1}{2}\delta_1^2$ and $\delta_1$ are cyclic 1-cocycles
on $\mathcal{H}_1$, and
$$F:=X\otimes Y-Y\otimes X -\delta_1Y\otimes Y$$
is a  cyclic 2-cocycle. See \cite{cm5} for detailed calculations and  relations
between these cocycles and
 the Schwarzian derivative, Godbillon-Vey cocycle, and
 the transverse fundamental class of Connes \cite{c3}, respectively.

\subsection{ Hopf-cyclic cohomology: type A, B, and C theories}

We recall the definitions of the three cyclic cohomology theories that were defined in \cite{hkrs2}.
We call them $A$, $B$ and $C$ theories. In the first case the algebra $A$ is endowed with
an action of a Hopf
algebra; in the second case the algebra $B$ is equipped with a coaction of a Hopf algebra; and finally
in theories of type $C$, we have a coalgebra endowed with  an action of a Hopf algebra. In all three
theories  the module of coefficients is a stable anti-Yetter-Drinfeld (SAYD)
 module over the Hopf algebra and we attach a
 cocyclic module in the sense of Connes \cite{c1} to the given
data.  Along the same lines
one can define a Hopf-cyclic cohomology theory for comodule coalgebras
as well (type D theory). Since so far we have found no applications of
such a theory we won't give its definition here. We also show that all
known examples of cyclic cohomology theories that are introduced so far such as: ordinary
cyclic cohomology for algebras, Connes-Moscovici's cyclic cohomology for
Hopf algebras \cite{cm1}, and equivariant cyclic cohomology \cite{ak1, ak2},  are  special cases  of
these theories.

Let $A$ be a left $H$-module algebra  and  $M$
be a left-right SAYD  $H$-module.
 Then the spaces $M\otimes A^{\otimes (n+1)}$ are right $H$-modules
via the diagonal action
$$(m\otimes \widetilde{a})h:=mh^{(1)}\otimes S (h^{(2)})\widetilde{a},$$
where the left  $H$-action on $\widetilde{a} \in A^{\otimes (n+1)}$ is
 via the left diagonal action of $H$.

We  define the space of {\it equivariant  cochains on $A$ with
coefficients in $M$} by
$$\mathcal{C}^n_H(A, M):= Hom_H(M\otimes A^{\otimes (n+1)}, k).$$
More explicitly, $f: M\otimes A^{\otimes (n+1)} \rightarrow k$ is in
$\mathcal{C}^n_H(A, M)$, if and only if
$$f((m\otimes a_0\otimes \cdots \otimes a_n)h)= \varepsilon (h) f(m\otimes a_0\otimes \cdots \otimes a_n), $$
for all $h \in H, m\in M$, and $a_i\in A$. It is shown in \cite{hkrs2}
that the following operators define a cocyclic module structure on $\{\mathcal{C}^n_H(A, M)\}_{n\in \mathbb{N}}$:
\begin{align*}&
(\delta_if)(m\otimes a_0\otimes \cdots\otimes a_n)
= f(m\otimes a_0\otimes \cdots \otimes a_i a_{i+1}\otimes\cdots \otimes a_n),
~~~ 0 \leq i <n,~~~~~~\\ &
(\delta_nf)(m\otimes a_0\otimes \cdots\otimes a_n)
=f(m^{(0)}\otimes (S^{-1}(m^{(-1)})a_n)a_0\otimes a_1\otimes\cdots \otimes a_{n-1}),~~~~~~
\\ &
(\sigma_if)(m\otimes a_0\otimes \cdots\otimes a_n)
= f(m\otimes a_0 \otimes \cdots \otimes a_i\otimes 1 \otimes \cdots \otimes a_n),
~~~0\le i\le n,~~~~~~
\\  &
(\tau_nf)(m\otimes a_0\otimes \cdots\otimes a_n)
= f(m^{(0)}\otimes S^{-1}(m^{(-1)})a_n\otimes a_0\otimes \cdots \otimes a_{n-1}).~~~~~~
\end{align*}

We denote the resulting cyclic cohomology theory by $HC^n_H(A, M), n=0,
1, \cdots.$

\begin{example}
 1. For $H=k=M$ we obviously recover the  standard cocyclic module of the  algebra $A$. The resulting
 cyclic cohomology theory
is the ordinary cyclic cohomology of algebras.\\

\noindent  2. For $M=H$ and $H$ acting on $M$ by conjugation and coacting
via  coproduct (Example 2.2.3.), we obtain the equivariant cyclic
cohomology theory of Akbarpour and Khalkhali For $H$-module algebras \cite{ak1, ak2}. \\

\noindent 3. For $H=k[\sigma, \sigma^{-1}]$ the Hopf algebra of Laurent polynomials, where $\sigma$ acts by
automorphisms on an algebra $A$, and $M=k$ is a trivial module, we obtain the so called
{\it twisted cyclic
cohomology} of $A$ with respect to $\sigma$. A {\it twisted cyclic n-cocycle} is a
linear map  $f: A^{\otimes (n+1)} \to k$  satisfying:
    $$f(\sigma a_n, a_0, \cdots,  a_{n-1})= (-1)^n f(a_0, \cdots,
    a_n), \quad  b_\sigma f=0,$$
    where $b_\sigma$ is the twisted Hochschild boundary defined by
\begin{eqnarray*}
b_\sigma f (a_0, \cdots, a_{n+1})&=&\sum_{i=0}^n (-1)^if(a_0, \cdots,
a_ia_{i+1}, \cdots, a_{n+1})\\
&+&(-1)^{n+1} f(\sigma (a_{n+1})a_0, a_1, \cdots,
 a_{n}).
\end{eqnarray*}

 \noindent 4. It is easy to see that for $M=^{\sigma}\!\!\!k_{\delta}$, the SAYD
 module attached to a modular pair in involution $(\delta, \sigma)$, $HC^0_H(A, M)$ is the space of $\delta$-invariant $\sigma$-traces
 on $A$ in the sense of Connes-Moscovici \cite{cm2,cm3} (cf. Section
 1.3.).
\end{example}

Next, let $B$ be a right $H$-comodule algebra and $M$ be a right-right SAYD $H$-module. Let
$${\mathcal C}^{n,H}(B,M):=Hom^H(B^{\otimes(n+1)},M),$$
denote the space of right $H$-colinear  $(n+1)$-linear functionals  on $B$ with values
in $M$. Here $B^{\otimes(n+1)}$ is considered a right $H$-comodule via
the
diagonal coaction  of $H$:
$$ b_0\otimes \cdots \otimes b_n \mapsto (b_0^{(0)}\otimes \cdots
\otimes b_n^{(0)})\otimes (b_0^{(1)} b_1^{(1)}\cdots b_n^{(1)}).$$
It is shown in \cite{hkrs2} that, thanks to the invariance
property imposed on our cochains and the SAYD condition on $M$,
the following maps define  a cocyclic module structure
 on $\{{\mathcal C}^{n,H}(B,M)\}_{n\in {\mathbb N}}$:
\begin{align*}
&(\delta_if)(b_0, \cdots, b_{n+1})
=f(b_0,  \cdots, b_i b_{i+1},\cdots, b_{n+1}),~~~
0\le i< n,\\
&(\delta_{n}f)(b_0, \cdots, b_{n+1})=
f(b_{n+1}^{(0)}b_0,  b_1, \cdots, b_{n})b_{n+1}^{(1)},\\
&(\sigma_if)(b_0, \cdots, b_{n-1})=
f(b_0, \cdots, b_i, 1, \cdots  b_{n-1}),~~~
0\leq i< n-1,\\
&(\tau_n f)(b_0, \cdots, b_{n})=
f(b_n^{(0)}, b_0, \cdots,  b_{n-1})b_n^{(1)}.
\end{align*}
We denote the resulting cyclic cohomology groups by $HC^{n, H}(B, M)$,
$n=0, 1, \cdots$.

\begin{example}
1. For $B=H$, equipped with comultiplication as coaction, and
$M=^{\sigma}\!\!\!k_{\delta}$ associated to a modular pair in
involution, we obtain the Hopf-cyclic cohomology of Hopf algebras
as defined in \cite{kr1}. This theory is different from
Connes-Moscovici's theory for Hopf algebras.
  It is dual, in the sense of Hopf algebras and not Hom dual,  to Connes-Moscovici's
  theory \cite{kr3}. It is computed in the
  following cases \cite{kr1}: $H=kG$, $H=U(\mathfrak{g})$, where it is isomorphic to group cohomology
  and Lie algebra
  cohomology, respectively;   $H=SL_2(q, k)$, and  $H=U_q(sl_2).$
  \\

\noindent 2. For $H=k$, and $M=k$ a trivial module, we obviously recover the cyclic cohomology of the
algebra $B$.

\end{example}

Finally we describe theories of type $C$ and their main examples. As we
shall see, Connes-Moscovici's original example of Hopf-cyclic cohomology
belong to this class of theories.
Let $C$ be a left $H$-module coalgebra, and $M$ be a right-left SAYD $H$-module. Let
$${\mathcal C}^n(C,M):=M\otimes_H C^{\otimes(n+1)} \quad n\in{\mathbb N}.$$
It can be checked that, thanks to the SAYD condition on $M$,   the
following operators are well defined and define a cocyclic module, denoted
$\{{\mathcal C}^n(C,M)\}_{n\in{\mathbb N}}$. In particular the
crucial periodicity conditions $\tau_n^{n+1}=id, \, n=0, 1, 2 \cdots$,
are  satisfied \cite{hkrs2}:
\begin{eqnarray*}
\delta_i(m\otimes c_0 \otimes \cdots \otimes c_{n-1})
&=&
m\otimes  c_0 \otimes\cdots\otimes  c_i^{(1)}\otimes c_i^{(2)}\otimes
c_{n-1},\,
 0 \leq i <n,\\
\delta_{n}(m\otimes c_0 \otimes \cdots \otimes c_{n-1})
&=&
m^{(0)}\otimes  c_0^{(2)}\otimes c_1
\otimes \cdots \otimes c_{n-1}  \otimes m^{({-1})}c_0^{(1)},
\\
\sigma_i(m\otimes c_0 \otimes \cdots \otimes c_{n+1})
&=&
m\otimes c_0 \otimes \cdots
\otimes \varepsilon(c_{i+1})\otimes
\cdots\otimes c_{n+1},\,
0\leq i \leq n,\\
\tau_n(m\otimes c_0 \otimes  \cdots \otimes c_n)
&=&
m^{(0)}\otimes c_1 \otimes
\cdots \otimes c_n \otimes m^{({-1})}c_0.
\end{eqnarray*}

\begin{example}
1. For $H=k=M$, we recover the cocyclic module of a coalgebra which defines its cyclic cohomology.\\

\noindent 2. For $C=H$ and $M=^{\sigma}\!\!\!k_{\delta}$, the
cocyclic module $\{{\mathcal C}^n_H(C,M)\}_{n\in{\mathbb N}}$ is
isomorphic to  the cocyclic module of Connes-Moscovici \cite{cm2},
attached to a Hopf algebra endowed with a modular pair in
involution. This example is truly fundamental and started the
whole theory.
\end{example}


\begin{thebibliography}{9}
\bibitem{ak1}
R. Akbarpour, and M. Khalkhali, { Hopf algebra equivariant cyclic homology and cyclic homology of
 crossed product algebras.} J. Reine Angew. Math., {\bf 559} (2003), 137-152.

\bibitem{ak2}
 R. Akbarpour, and M. Khalkhali, Equivariant cyclic cohomology of $ H$-algebras.
 $K$-Theory  {\bf 29}  (2003),  no. 4, 231--252.

\bibitem{c1}
A. Connes, {Cohomologie cyclique et foncteurs  ${\rm Ext}\sp n$ }. (French) (Cyclic cohomology
and functors ${\rm Ext}\sp n$) C. R.
 Acad. Sci. Paris Sér. I Math. {\bf 296} (1983), no. 23, 953--958.

 \bibitem{c2}
 A. Connes, {Noncommutative differential geometry}. Inst. Hautes Études Sci.
 Publ. Math. No. {\bf 62} (1985), 257--360.

\bibitem{c3}
 A. Connes,  Cyclic cohomology and the transverse fundamental class of a foliation.
 Geometric methods in operator algebras (Kyoto, 1983),  52--144, Pitman Res. Notes Math. Ser., 123,
 Longman Sci. Tech., Harlow, 198.

\bibitem{cma1} A. Connes, and M. Marcolli, From physics to number theory via noncommutative
geometry. Part I: Quantum statistical mechanics of Q-lattices.
math.NT/0404128.

\bibitem{cma2} A. Connes, and M. Marcolli,
 From physics to number theory via noncommutative geometry, Part II: Renormalization,
 the Riemann-Hilbert correspondence, and motivic Galois theory.
 hep-th/0411114.


\bibitem{cm1} A. Connes and H. Moscovici, {Local index formula in
noncommutative geometry}, Geom. Funct. Anal. {\bf 5} (1995), no. 2,
174-243.

\bibitem{cm2} A. Connes and H. Moscovici, {Hopf algebras, Cyclic Cohomology and the transverse
index theorem}, Comm. Math. Phys. {\bf 198} (1998), no. 1, 199--246.

\bibitem{cm3}
   A. Connes and H. Moscovici, {Cyclic cohomology and Hopf algebras}.
    Lett. Math. Phys.  {\bf 48}  (1999),  no. 1, 97--108.


\bibitem{cm4}
A. Connes, and   H. Moscovici, {Cyclic cohomology and Hopf algebra symmetry.}
Conference Moshé Flato 1999 (Dijon).
 Lett.
Math. Phys. {\bf 52} (2000), no. 1, 1--28.

\bibitem{cm5}
A. Connes and   H. Moscovici, Modular Hecke algebras and their Hopf
symmetry,
 to appear in the Moscow Mathematical Journal (volume dedicated to
 Pierre Cartier).

\bibitem{cm6}
A. Connes and   H. Moscovici, Rankin-Cohen brackets and the Hopf algebra
of transverse geometry,
to appear in the Moscow Mathematical Journal (volume dedicated to Pierre
Cartier).




 \bibitem{hkrs1}
P. M. Hajac, M. Khalkhali, B. Rangipour, and Y. Sommerh\"auser,
{Stable anti-Yetter-Drinfeld modules.}  C. R. Math. Acad. Sci. Paris
{\bf 338}  (2004),  no. 8, 587--590.

\bibitem{hkrs2}
P. M. Hajac, M. Khalkhali, B. Rangipour, and Y. Sommerh\"auser,
{Hopf-cyclic homology and cohomology with coefficients.}  C. R. Math. Acad. Sci. Paris  {\bf 338}  (2004),
no. 9, 667--672.

\bibitem{js}
P. Jara, D. Stefan, Cyclic homology of Hopf Galois extensions and Hopf algebras.
Preprint math/0307099.



\bibitem{js}
A. Joyal, and R. Street, The geometry of tensor calculus. I.  Adv.
Math.  {\bf 88}  (1991),  no. 1, 55--112.

\bibitem{kas}
Ch. kassel,
Quantum groups. Graduate Texts in Mathematics, 155. Springer-Verlag, New
York, 1995.

\bibitem{kay}
A. Kaygun,  Bialgebra cyclic homology with coefficients. Preprint
math/0408094.


\bibitem{kr1}
M. Khalkhali, and  B. Rangipour, {A new cyclic module for Hopf algebras}.
$K$-Theory {\bf 27} (2) (2002), 111-131.

 \bibitem{kr2}
 M. Khalkhali, and  B. Rangipour, {Invariant cyclic homology}. $K$-Theory {\bf 28} (2) (2003), 183-205.
 

\bibitem{kr3}
 M. Khalkhali, and  B. Rangipour, {A note on cyclic duality and Hopf algebras}.
 Comm. in Algebra {\bf 33} (2005) 763-773.


 \bibitem{kr4}
 M. Khalkhali, and  B. Rangipour, {Cup products in Hopf-cyclic cohomology}.
    C.R. Acad. Sci. Paris Ser. I {\bf 340} (2005), 9-14.
 

\bibitem{mar}
M. Marcolli, Lectures on noncommutative arithmetic geometry.
math.QA/0409520.

\bibitem{ml}
 S. Mac Lane, categories for the working mathematician. Second edition.
 Graduate Texts in Mathematics, 5. Springer-Verlag, New York, 1998.

\bibitem{maj}
S. Majid, {Foundations of quantum group theory.} Cambridge University Press, Cambridge, 1995.

\bibitem{rt}
D. Radford, and J. Towber,  Yetter-Drinfeld categories associated to an arbitrary
bialgebra. J. Pure Appl. Algebra {\bf 87} (1993), 259-279.

\bibitem{sta}
M. Staic, personal communication.

\bibitem{sw}
M. Sweedler, {Hopf algebras.} Mathematics Lecture Note Series W. A.
Benjamin, Inc., New York 1969.

\bibitem{ta}
 R. Taillefer, {Cyclic Homology of Hopf Algebras}. $K$-Theory {\bf 24} (2001), no. 1, 69--85.

\bibitem{voi}
C. Voigt, personal communication.

\bibitem{y}
D. Yetter,  Quantum groups and representations of monoidal categories.
Math. Proc. Cambridge Philos. Soc.  {\bf 108}  (1990),  no. 2, 261--290.
\end{thebibliography}
\end{document}